\documentclass[11pt,a4paper]{article}

\usepackage{bezier,amsfonts,amssymb,graphicx,amsthm,amsmath}
\usepackage[english]{babel}
\def\lesta{ \hfill $\Box$ \bigskip}

\newcommand{\auttfg}{$\mbox{Aut}^{\mbox{{\tiny{\textbf{TF}}}}}(G)$}
\newcommand{\auttfgg}{\emph{Aut}$^{\mbox{\tiny\textbf{{TF}}}}(G)$}

\newcommand{\BG}{\mbox{\textbf{CDC}}(G)}

\newcommand{\IDC}{\mbox{\textbf{IDC}}}
\newcommand{\ADC}{\mbox{\textbf{ADC}}}

\newcommand{\ZZ}{\mathbb Z}

\newtheorem{Thm}{Theorem}[section]

\newtheorem{Prop}[Thm]{Proposition}
\newtheorem{Cor}[Thm]{Corollary}
\newtheorem{Lem}[Thm]{Lemma}

\begin{document}
\title{A Generalisation of Isomorphisms with Applications }

\author{ J. Lauri, R. Mizzi \\ Department of Mathematics \\ University of Malta
\\ Malta \\ josef.lauri@um.edu.mt \\ russell.mizzi@um.edu.mt \and R. Scapellato  \\ Dipartimento 
di Matematica \\ Politecnico di Milano \\ Milano \\ Italy \\ raffaele.scapellato@polimi.it}
\pdfinfo{ /Author (Josef Lauri, Russell Mizzi, Raffalele Scapellato) /Title(Two-fold Orbital mixed graphs) /Keywords (graph theory) }\maketitle

\begin{small}
\begin{abstract}
In this paper, we study the behaviour of TF-isomorphisms, a natural generalisation of isomorphisms. TF-isomorphisms allow us to simplify the approach to seemingly unrelated problems. In particular, we mention the Neighbourhood Reconstruction problem, the Matrix Symmetrization problem and Stability of Graphs. We start with a study of invariance under TF-isomorphisms. In particular, we show that alternating trails and incidence double covers are conserved by TF-isomorphisms, irrespective of whether they are TF-isomorphisms between graphs or digraphs. We then define an equivalence relation and subsequently relate its equivalence classes to the incidence double cover of a graph. By directing the edges of an incidence double cover from one colour class to the other and discarding isolated vertices we obtain an invariant under TF-isomorphisms which gathers a number of invariants. This can be used to study TF-orbitals, an analogous generalisation of the orbitals of a permutation group. \end{abstract}
\end{small}
 \bigskip


\section{Introduction}

Consider the two graphs shown in Figure \ref{fig:peterporcu}. One is the well-known Petersen graph which we denote by $\Pi$  and the other is a graph which is not so well-known which we sometimes refer to as Petersen's cousin, for reasons which will soon become apparent, and which we denote by $\Lambda$. What relationship could there be between them? Consider the set of neighbourhoods of the vertices in the two graphs. These are:

\begin{figure}[h]
\begin{centering}
 \includegraphics[width=9.3cm,height=4.25cm]{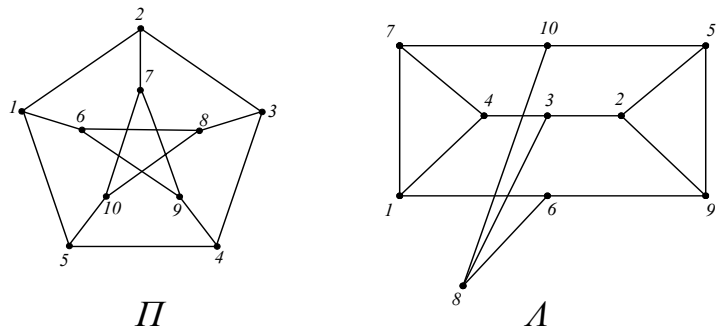}
 \caption{The Petersen graph $\Pi$ and its less well-known cousin $\Lambda$.} \label{fig:peterporcu}
\end{centering} 
\end{figure}

\bigskip\noindent
{Neighbourhoods of $\Pi$: \{2,5,6\},\ \{1,3,7\},\  \{2,4,8\},\  \{3,5,9\},\  \{1,4,10\},\  \{1,8,9\},\  \{2,9,10\},\  \{3,6,10\},\  \{4,6,7\},\  \{5,7,8\}. }\\
{Neighbourhoods of $\Lambda$: \{4,6,7\},\  \{3,5,9\},\  \{2,4,8\},\  \{1,3,7\},\  \{2,9,10\},\  \{1,8,9\},\  \{1,4,10\},\  \{3,6,10\},\  \{2,5,6\},\  \{5,7,8\}. }\\

\bigskip
Up to a re-ordering, both graphs have the same family of neighbourhoods. It is therefore clear that if one were given just the family of neighbourhoods of the Petersen graph one would not be able to determine that the graph they came from was Petersen---it could have been the second graph which also has the same neighbours. \\

\bigskip

In the literature the following problem (the Neighbourhood Reconstruction Problem) has been proposed (for example, in  \cite{Aigner1} and \cite{Aigner2}): given the neighbourhoods of the vertices of of $G$, can $G$ be determined uniquely up to isomorphism? The two graphs above clearly show that the answer to this question is ``no'' in general. The Petersen graph is not reconstructible this way because the second graph shown in the figure is a reconstruction of the Petersen which is not isomorphic to it. Why does this happen? We shall explain this below. How many other reconstructions of the Petersen graph can be obtained this way? We shall see later on that this second graph is the only such reconstruction of the Petersen graph. \\

There are a few other problems which have been considered in the graph theory literature which, as we shall see, are closely related to the neighbourhood reconstruction problem.\\

\begin{enumerate}
\item {\bf The Realisability Problem.} When is a given family of vertices the neighbourhood family of a graph or a digraph? What is the computational complexity of determining whether such a given family is the neighbourhood family of a graph or a digraph?\\

\item {\bf The Matrix Symmetrization Problem.} Given a $(0,1)$-matrix $A$, is it possible to change it into a symmetric matrix using (independent) row and column permutations? Although it is not immediately obvious, we shall see that this problem is related to the Realisability Problem. This problem was first studied in the paper \cite{Scapsalvi1} starting with a matrix $A$ which is already symmetric.\\

\item {\bf Stability.} This problem was first raised and studied in \cite{Scapsalvi1}. Given the categorical product $G \times K_2$ of a graph or digraph $G$ with the complete graph $K_2$, the graph $G$ is said to be \emph{unstable} when the automorphism group of the product is not isomorphic to Aut$(G) \times \ZZ$. When is a graph unstable?  This question was heavily studied in \cite{Scapsalvi2,Surowski1,Surowski2,Wilson01}, and again, although not immediately clear why, it is strongly related to the previous questions. \\

\end{enumerate}

\smallskip
The excellent survey paper \cite{Gurvich1} gives a good historical picture of work done on these problems.\\

\bigskip\noindent
In this paper we shall present a new type of isomorphism between graphs and digraphs which, we believe, has independent interest but also unifies the above problems, as we shall demonstrate along the way while presenting our results.\\


\section{Notation}

A \textit{mixed graph} is a pair $G=(\mbox{$V$}(G),\mbox{$A$}(G))$ where $V$$(G)$ is a set and $A(G)$ is a set of ordered pairs of elements of $V$$(G)$. The elements of $V$$(G)$ are called \textit{vertices} and the elements of $A(G)$ are called \textit{arcs}. When referring to an arc $(u,v)$, we say that $u$ is \textit{adjacent to} $v$ and $v$ is \textit{adjacent from} $u$. The vertex $u$ is the \textit{tail} and $v$ is the \textit{head} of a given arc $(u,v)$. An arc of the form $(u,u)$ is called a \textit{loop}.  A mixed graph cannot contain multiple arcs, that is, it cannot contain the arc $(u,v)$ more than once.  A set $S$ of arcs is \textit{self-paired} if, whenever $(u,v) \in$ $S$, $(v,u)$ is also in $S$. If $S$ $\ =\{(u,v), (v,u)\}$, then we identify $S$ with the unordered pair $\{u,v\}$; this unordered pair is called an \textit{edge}.\\

It is useful to consider two special cases of mixed graphs. A \textit{graph} is a mixed graph without loops whose arc-set is self-paired. The edge set of a graph is denoted by $E$$
(G)$. A \textit{digraph} is a mixed graph with no loops in which no set of arcs is self-paired. The \textit{inverse} $G'$ of a mixed graph $G$ is obtained from $G$ by reversing all its arcs, that is $V$$(G') =$$V$$(G)$ and $(v,u)$ is an arc of $G'$ if and only if $(u,v)$ is an arc of $G$. A digraph $G$ may therefore be characterised as a mixed graph for which $A(G)$ and $A(G')$ are disjoint and a graph as one for which $A(G)=A(G')$. The underlying graph $\widehat{G}$ of a mixed graph $G$ is a graph with the vertex set $V$$(\widehat{G})$ $=$ $V$$(G)$ and the edge set $E$$(\widehat{G})$ defined by $\{x,y\} \in$ $E$$(\widehat{G})$ if and only if either $(x,y)$ or $(y,x)$ is an element of $A(G)$. Two arcs are \emph{incident} in $G$ if the corresponding edges in the underlying graph $\widehat{G}$ have a common vertex. When we say that a \emph{mixed graph is connected}, we mean that the underlying graph is connected.\\

Given a mixed graph $G$ and a vertex $v \in$ $V$$(G)$, we define the \textit{in-neighbourhood} $N_{in}(v)$ by $N_{in}(v) = \{x \in \mbox{$V$}(G)- (x,v) \in \mbox{A}(G)\}$. Similarly we define  the  \textit{out-neighbourhood} $N_{out}(v)$ by $N_{out}(v) = \{x \in \mbox{$V$}(G)-(v,x) \in \mbox{A}(G)\}$.  The \textit{in-degree} $ \rho_{in}(v)$ of a vertex $v$ is defined by $ \rho_{in}(v) = |N_{in}(v)|$ and the \textit{out-degree} $ \rho_{out}(v)$ of a vertex $v$ is defined by $ \rho_{out}(v) = |N_{out}(v)|$. When $G$ is a graph, these notions reduce to the usual neighbourhood $N(v)=N_{in}(v)=N_{out}(v)$ and degree $\rho(v)=\rho_{in}(v)=\rho_{out}(v)$. A vertex $v$ is called a \textit{source} if  $ \rho_{in}(v)= 0$ and a \textit{sink} if $ \rho_{out}(v)=0$. A vertex is said to be \textit{isolated} when it is both a source and a sink, that is, it is not adjacent to or from any vertex.\\

A mixed graph $G$ is called \textit{bipartite} if there is a partition of $V$$(G)$ into two sets $X$ and $Y$, which we call \textit{colour classes}, such that for each arc $(u,v)$ of $G$ the set $\{u,v\}$ intersects both $X$ and $Y$.  We call a bipartite digraph having one colour class consisting of sources and the other colour class consisting of sinks as a \emph{strongly bipartite digraph}. \\

Let $G$ be a digraph and let $(u,v)$ be an arc of $G$. If in $G-(u,v)$, the vertices $u$, $v$ are either both sources or both sinks, then we call $(u,v)$ an S-\textit{arc} of $G$.\\

A set $P$ of arcs of $G$ is called a \emph{trail} if its elements can be ordered in a sequence $a_{1},\ a_{2}, \dots,\ a_{k}$ such that each $a_{i}$ is incident with $a_{i+1}$ for all $i = 1,\ \dots,\ k-1$. If $u$ is the vertex of $a_{1}$, that is not in $a_{2}$ and $v$ is the vertex of $a_{k}$ which is not in $a_{k-1}$, then we say that $P$ \emph{joins} $u$ and $v$; $u$ is called the \emph{first vertex} of $P$ and $v$ is called the \emph{last vertex} with respect to the sequence $a_{1},\ a_{2},\ \dots,\ a_{k}$. If, whenever $a_{i}=(x,y)$, either $a_{i+1}=(x,z)$ or  $a_{i+1}=(z,y)$ for some new vertex $z$, $P$ is called  an \textit{alternating trail} or \textbf{A}-\textit{trail}.\nocite{Zelinka2}  \\

If the first vertex $u$ and the start-vertex $v$ of an \textbf{A}-trail $P$ are different, then $P$ is said to be \emph{open}. If they are equal then we have to distinguish between two cases. When the number of arcs is even then $P$ is called \emph{closed} while when the number of arcs is odd then $P$ is called \emph{semi-closed}. Note that if $P$ is semi-closed then either (i) $a_{1}=(u,x)$ for some vertex $x$ and $a_{k} = (y,u)$ for some vertex $y$ or (ii) $a_{1}=(x,u)$ and $a_{k} = (u,y)$. If $P$ is closed then either $a_{1} =(u,x)$ or $a_{k}=(u,y)$ or $a_{1}=(x,u)$ and $a_{k} = (y,u)$.  Observe also that the choice of the first (equal to the last) vertex for a closed \textbf{A}-trail is not unique but depends on the ordering of the arcs. However, this choice is unique for semi-closed \textbf{A}-trails as this simple argument shows. Suppose $P$ is semi-closed and the arcs of $P$ are ordered such that $u$ is the unique (in that ordering) first and last vertex, that is, it is the unique vertex such as the first and the last arcs in the ordering in $P$ do not alternate in direction at the meeting point $u$. Therefore, it is easy to see that both $\rho_{in}(u)$ and $\rho_{out}(u)$ (degrees taken in $P$ as a subgraph induced by its arcs) are odd whereas any other vertex $v$ in the trail has both $\rho_{in}(v)$ and $\rho_{out}(v)$ even. This is because, in the given ordering, arcs have to alternate in direction at $v$ and therefore in-arcs of the form $(x,v)$ are paired with out-arcs of the form $(v,y)$. Therefore, in no ordering of the arcs of $P$ can $u$ be anything but the only vertex at which the first and last arcs do not alternate. The same argument holds for open \textbf{A}-trails. Therefore, open and semi-closed \textbf{A}-trails are similar at least in the sense that the first and last vertices are uniquely determined regardless of the sequence of the arcs. This similarity will be strengthened by the results which we shall shortly present.\\

Let $G$ be a mixed graph. If, for every $u$, $v \in$ $V$$(G)$, there exists an \textbf{A}-trail which joins them, then we say that $G$ is \textbf{A}-\textit{connected}.   Clearly, a connected graph $G$ with at least two edges is always \textbf{A}-connected since we can always choose any orientation for a given edge. In the case of a mixed graph, \textbf{A}-connectedness is not guaranteed. In fact, there are easy counterexamples also among digraphs. \\

In a connected graph, the length (that is, number of edges) of a shortest path between two given vertices $u$, $v$ is denoted by $d(u,v)$. Any other graph theoretical terms which we use  are standard and can be found in textbooks such as \cite{bondy} and \cite{Harary01}. Information on automorphism groups of a graph can be found in \cite{lauri2}.\\

Let $G$ and $H$ be two mixed graphs and $\alpha$, $\beta$ be bijections from $V$$(G)$ to $V$$(H)$. The pair $(\alpha,\beta)$ is said to be a \textit{two-fold 
isomorphism} (or TF\textit{-isomorphism}) if the following holds: $(u,v)$ is an arc of $G$ if and only if $(\alpha(u),\beta(v))$ is an arc of $H$. We then say that $G$ and $H$ are  TF-\textit{isomorphic} and write $G\cong ^{\mbox{{\tiny{\textbf{TF}}}}} H$. Note that when $\alpha=\beta$ the pair $(\alpha,\beta)$ is a TF-isomorphism if and only if $\alpha$ itself is an isomorphism. If $\alpha \neq \beta$, then the given TF-isomophism $(\alpha,\beta)$ is essentially different from a usual isomorphism and hence we call $(\alpha,\beta)$ a \textit{non}-\textit{trivial} TF-\textit{isomorphism}.  In this case, we also say that $G$ and $H$ are \emph{non-trivially} TF-\emph{isomorphic}. If $(\alpha,\beta)$ is a non-trivial TF-isomorphism from a mixed graph $G$ to a mixed graph $H$, the bijections $\alpha$ and $\beta$ need not necessarily be isomorphisms from $G$ to $H$. This is illustrated by examples found in \cite{lms2}, and also others found below. \\
When $G=H$, $(\alpha,\beta)$ is said to be a TF-\textit{automorphism} and it is again called non-trivial if $\alpha \neq \beta$. The set of all TF-automorphisms of $G$ with 
multiplication defined by $(\alpha,\beta)(\gamma,\delta) = (\alpha \gamma, \beta \delta)$ is a subgroup of $S_{V(G)} \times S_{V(G) }$ and it is called the \textit{two-fold 
automorphism group} of $G$ and is denoted by \auttfg.  Note that if we identify an automorphism $\alpha$ with the TF-automorphism $(\alpha,\alpha)$, then Aut$(G) \subseteq$ 
\auttfg. When a graph has no non-trivial TF-automorphisms, Aut$(G)= $\auttfg.  It is possible for an asymmetric graph $G$, that is a graph with $|$Aut$(G)| = 1$, to have non-trivial TF-automorphisms \cite{lms2}.\\

\section{Some double covers and invariants under TF-isomorphisms}

Let $G$ be a mixed graph. The \emph{incidence double cover} of $G$, denoted by \IDC$(G)$ is a bipartite graph with vertex set $V$(\IDC$(G)$)$\subseteq$ $V$$(G) \times \{0,\ 1\}$ and edge set $E$(\IDC$(G)$) $=$ $\{(u,0),(v,1)\} \ | \ (u,v) \in$ $A(G)\}$. The reader may refer to \cite{klin112} for more information regarding the incidence double cover of graphs and its relevance to the study of association schemes. The \textbf{A}-\emph{cover} of $G$, denoted by $\mbox{\textbf{ADC}}(G)$ is a strongly bipartite digraph with vertex set $V$$(\mbox{\textbf{ADC}}(G))$ $\subseteq$ $V$$(G) \times \{0,\ 1\}$ and arc set $A$($\mbox{\textbf{ADC}}(G)$) $=$ $\{(u,0),(v,1)\} \ $ $| \ (u,v)$ $ \in$ $A(G)\}$. For a more concise notation, very often we use $u_{0}$ or $u_{1}$ to label the elements of $V$$(G) \times \{0,\ 1\}$ instead of $(u,0)$ or $(u,1)$. It is clear that \textbf{ADC}$(G)$ is obtained from \textbf{IDC}$(G)$ by removing isolated vertices and changing every edge $\{u_{0},v_{1}\}$ into an arc $(u_{0},v_{1})$.\\

{\Thm{Let $G$, $H$ be mixed graphs. Then $G$ and $H$ are two-fold isomorphic if and only if \emph{\IDC$(G)$} and \emph{\IDC$(H)$} are isomorphic.}\label{thm:idctf}}\\

{\proof{Let $(\alpha,\beta)$ be any two-fold isomorphism from $G$ to $H$.
This implies that, for any $(u,v) \in$ $A(G)$, $(\alpha(u),\beta(v)) \in$ $A(H)$.
Consequently, given the corresponding $\{(u,0),(v,1)\} \in$ $E$(\IDC$(G))$, $\{(\alpha(u),0),(\beta(v),1)\}$ $\in$ $E$(\IDC$(H))$.
Define $\phi:$ $V$(\IDC$(G)$) $\rightarrow$ $V$(\IDC$(H)$) such that $\phi(x,0) = (\alpha(x),0)$ and $\phi(x,1) = (\beta(x),1)$
for any $x \in$ $A(G)$ such that $\phi$ is an isomorphism from \IDC$(G)$ to \IDC$(H)$.\\

Conversely, let $\alpha$ be any isomorphism from \IDC$(G)$ to \IDC$(H)$ such that $\alpha \{(u,0),$ $(v,1)\}$ $=$ $\{(u',0),$ $(v',1)\}$.
Clearly $(u,v)\in$ $A(G)$ and $(u',v') \in$ $A(H)$ by the definition of IDC. Define $\alpha :$ $V$($G$) $\rightarrow$ $V$($H$)  such that
$\alpha(u) = u'$ if and only if $\phi(u,0) = (u',0)$. Similarly define $\beta :$ $V$($G$) $\rightarrow$ $V$($H$) such that 
$\beta(v) = v'$ if and only if $\phi(v,1) = (v',1)$. Given any $(x,y) \in$ $A(G)$, $(\alpha(x),\beta(y)) \in$ $A(H)$ since $\{(x,0),(y,1)\} \in$ $E$(\IDC$(G)$) if and only if $\{\phi(x,0),\phi(y,1)\} = \{(x',0),(y',1)\}$ in $E$(\IDC$(H)$) if and only if $(x',y')$ $\in$ $A(H)$.
 }\lesta}

We now present what was Theorem 3.7 in \cite{lms1}, one of our main results in \cite{lms1},  as a corollary to Theorem \ref{thm:idctf}.\\

The \emph{canonical double cover (CDC)} of a graph or digraph $G$ (also called its \emph{duplex} especially in computational chemistry literature, for example, \cite{randic83}) is the graph or digraph whose vertex set is $V(G) \times \{0,1\}$ and in which there is an arc from $(u,i)$ to $(v,j)$ if and only if  $i \neq j$ and there is an arc from $u$ to $v$ in $G$. The canonical double cover of $G$ is often described as the direct or categorical product $G \times K_{2}$ \cite{imrich,HB},  and is sometimes also called the \emph{bipartite double cover} of $G$. For graphs, the canonical double cover is identical to the incidence double cover.\\

{\Cor{Two graphs $G$, $H$ are TF-isomorphic if and only if  \emph{\textbf{CDC}}$(G)$ and  \emph{\textbf{CDC}}$(H)$ are isomorphic.}\label{cor:tfidc01}}\\

{\proof{In fact, since $G$ and $H$ are graphs, \IDC$(G)$ $\cong$  \textbf{CDC}$(G)$ and \IDC$(H)$ $\cong$  \textbf{CDC}$(H)$. 
}\lesta}

Therefore, in general the IDC of a mixed graph $G$ is a structure which is invariant under the action of a TF-isomorphism acting on $G$. In the case of mixed graphs which are not graphs, Theorem \ref{thm:idctf} is a significant improvement over Theorem3.7 in \cite{lms1} which only considered TF-isomorphic graphs (not mixed graphs) and the canonical double cover. \\

 \bigskip
 
\subsection{Digression }

Now we can see how the neighbourhood reconstruction problem and the other problems we discussed in the first section can be described in terms of TF-isomorphisms. First, consider this alternative way of looking at TF-isomorphisms. An \emph{incidence structure} or, alternatively, a \emph{hypergraph}, is a finite set of vertices with a system of subsets (blocks) some of which can be repeated. Number the $n$ vertices of a hypergraph in some arbitrary but fixed way, and do similarly for the $b$ blocks of the hypergraph. The \emph{incidence matrix} of the hypergraph is the $n\times b$ matrix whose $ij$ entry is 1 if the $i$th vertex is in the $j$ block, and is zero otherwise. Let $H_1, H_2$ be two such hypergraphs with incidence matrices $B_1, B_2$, respectively. Then usually $H_1$ and $H_2$ are said to be isomorphic if there is a bijection $\alpha$ from $V(H_1)$ to $V(H_2)$ (effectively, a relabelling of the vertices of $H_1$) such that, under the resulting relabelling, the blocks of $\alpha(H_1)$ are the same as the blocks of $H_2$, possibly in a different order. Similarly, an automorphism of a hypergraph $H$ is a permutation of $V(H)$ (a relabelling of the vertices) such that the new blocks are a re-odering of the old blocks.\\

In other words, we have a permutation $\alpha$ of the rows of the incidence matrix $B_1$ such that the columns become a permutation of the columns of $B_2$. We can remove this last detail and make even the columns the same as those of $B_2$ by saying that an isomorphism from $H_1$ to $H_2$ is an independent re-ordering $\alpha$ of the rows and $\beta$ of the columns of $B_1$ such that it becomes $B_2$. Similarly, an automorphism of $H$ is an independent re-ordering of the rows and columns of $B$ which leaves $B$ unchanged. Therefore if we consider the adjacency matrix $A$ of a graph $G$ as an incidence matrix of a hypergraph with $n$ vertices (corresponding to the rows) and $b=n$ blocks (corresponding to the columns), a TF-isomorphism (TF-automorphism) is an isomorphism (automorphism) of the hypergraph represented by $A$.\\

Looking back at the example of the Petersen graph $\Pi$ and what we have called its cousin $\Lambda$ we see that their neighbourhoods considered as the blocks of two hypergraphs give isomorphic hypergraphs which means, according to the previous discussions, that $\Pi$ and $\Lambda$ are non-trivially TF-isomorphic, and that is why one is a neighbourhood reconstruction of the other! What non-trivial TF-isomorphism can we write from one to the other?   Looking at how the list of neighbourhoods of the vertices $\{1,2,\ldots,10\}$ of the second graph appear as a permutation of the same list of neighbourhoods of the first graph easily indicates that if $\alpha=\mbox{id}$ and $\beta=(1\ 9)(2\ 4)(5\ 7)$ then $(\alpha,\beta)$ is a TF-isomorphism from the $\Pi$ to $\Lambda$ as labelled in Figure \ref{fig:peterporcu}.\\

But how do we know that $\Lambda$ is the only graph which is a neighbourhood reconstruction of (that is, TF-isomorphic to) the Petersen graph? We shall soon see this when below, we present one more result on canonical double covers.  \\

The Matrix Symmetrization Problem can also be described in terms of TF-isomorphisms: given a digraph $D$, is there a graph  $G$ to which $D$ is non-trivially TF-isomorphic? In the case when the matrix $A$ is already symmetric, as the problem was originally posed in \cite{Scapsalvi1}, this question becomes: given a graph $G$ is it non-trivially isomorphic to some other graph (possibly $G$ itself)?\\

\begin{figure}
\begin{centering}
 \includegraphics[width=8cm,height=7.5cm]{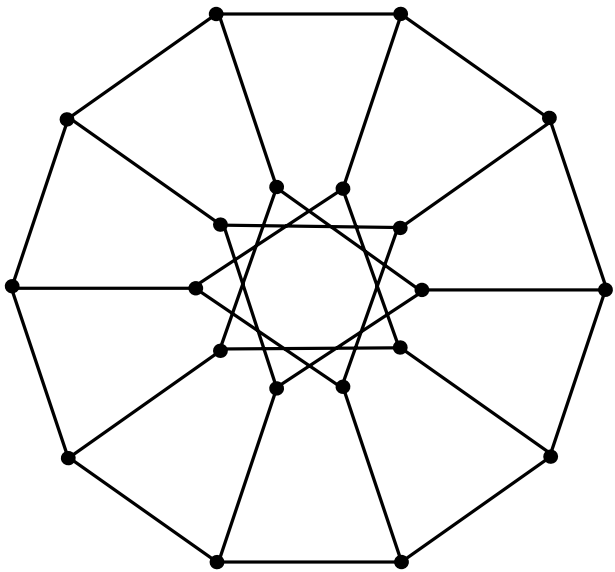}
 \caption{The Desargues graph.} \label{fig:desargues}
\end{centering} 
\end{figure}

Let us now return to the Petersen graph $\Pi$ and its cousin $\Lambda$  from Figure \ref{fig:peterporcu}. Since these two graphs are TF-isomorphic then they have the same CDC, and in fact, their common CDC is the well known Desargues graph shown in Figure \ref{fig:desargues}. (We have labelled Petersen's cousin by $\Lambda$ in honour of Livio Porcu who seems to have been the first one to observe in \cite{porcu} that $\Pi$ and $\Lambda$ have the same CDC.)\\

But now we can explain why these two graphs are the only ones with the same neighbourhood family. First we recall this result proved by Pacco and Scapellato in \cite{pacco}. As an easy reference, Theorem 5.3 of  \cite{pacco} may be restated as follows using our current terminology.\\

{\Thm{Given a connected bipartite mixed graph $H$, the number of non-isomorphic mixed graphs $G$ such that ${\emph{\BG}} \cong H$  is equal to the number of conjugacy classes of involutions in \emph{Aut}$(H)$ that interchange the two colour classes of $H$. The number  of non-isomorphic loopless mixed graphs $H$ such that ${\emph{\BG}} \cong H$ is equal to the number of conjugacy classes of involutions in \emph{Aut}$(H)$ that interchange the two colour classes of   $H$ \emph{and} do not take any vertex $u$ to a vertex $v$ such that $(u,v)$ is an arc.}\label{thm:conjugacyclassespeterporcu}\lesta}

Now, the automorphism group of the Desergues graph $D$ is isomorphic to $S_5 \times Z_2$, and has order 240. Letting $\beta$ be the automorphism of $D$ taking $(v,0)$ into $(v,1)$ and vice versa, note that $\beta$ belongs to the centre of the group. Hence, each involution of Aut$(D)$ takes the form  $(\alpha,{\rm id})$ or $(\alpha,\beta)$, where $\alpha$ is an involution of $S_5$. Only the latter swaps the two colour classes; its conjugacy classes are as many as those of involutions of $S_5$. The number of conjugacy classes of involutions of $S_5$ is exactly $2$, corresponding to transpositions and double transpositions. Therefore, by Theorem \ref{thm:conjugacyclassespeterporcu}, there are exactly two non-isomorphic mixed graphs whose CDC is $D$. One of them must be Petersen itself, while the other one is $\Lambda$, for which we know already that the CDC is $D$. So in this case, only these proper graphs occur, not more general mixed graphs.\\

Observe that since the Petersen graph's automorphism group is isomorphic to $S_5$, this graph is stable. However, Aut($\Lambda$) is isomorphic to
$S_3\times Z_2$. Thus the index of the automorphism group of $\Lambda$ in Aut$(D)$ is $20$ and so it is unstable.\\



\section{TF-isomorphism and alternating trails}

We shall consider isomorphisms and TF-isomorphisms between pairs of mixed graphs, that is, we shall allow loops, directed arcs and edges. 
Configurations conserved by TF-isomorphisms must also be conserved by isomorphisms since the latter are just a special case of the former. However, the converse does not necessarily hold.   It is well known that loops, paths and cycles are all conserved by isomorphisms, but it is easy to see that they are not necessarily conserved by TF-isomorphisms. For example, an arc $(u,v)$ can be mapped into a loop by a TF-isomorphism $(\alpha,\beta)$ if $\alpha(u)=\beta(u)$.\\

An isomorphism conserves degrees, in-degrees and out-degrees. In the case of TF-isomorphisms, the situation is slightly more elaborate. First note that $\alpha$ must conserve  the out-degree of each vertex but not the in-degree. Likewise, $\beta$ must conserve the in-degree of each vertex but not the out-degree. Therefore, if some vertex $u \in$ $V(G)$ is a source, $\alpha(u)$ might not be a source but it is certainly not a sink. If $u \in$ $V(G)$ is a sink, then $\alpha(u)$ must also be a sink. An analogous argument  holds for $\beta$.  Hence in the case of a digraph whose vertex set consists only of sources and sinks, $\alpha$ and $\beta$ must take sources to sources and sinks to sinks.  Also, if $G$ and $H$ are graphs and  $(\alpha,\beta)$ is a TF-isomorphism from $G$ to some graph $H$, then $\alpha$, $\beta$ must preserve the degree $\rho(v)$ of any vertex $v \in$ $V(G)$ since for every vertex $v$ in a graph, $\rho(v) = \rho_{in}(v) = \rho_{out}(v)$.\\  

The definitions of the term \textit{path} found in the literature tacitly imply a specific direction from one  vertex to the subsequent vertex in a sequence. 
For example, if $u ,v, w$ is a path in graph $G$ and $\alpha$ is an isomorphism from $G$ to $H$, then the arcs $(u,v)$ and $(v,w)$ are mapped into the arcs $(\alpha(u),\alpha(v))$ and $(\alpha(v),\alpha(w))$ in $H$, with the common vertex $\alpha(v)$. But, if $(\alpha,\beta)$ is a TF-isomorphism from $G$ to $H$ then the arc $(u,v)$ is mapped into $(\alpha(u),\beta(v))$ and it is the arc $(v,w)$ which is mapped into $(\alpha(w),\beta(v))$ containing the common vertex $\beta(v)$ with the previous arc. That is, to obtain a common vertex between images of successive arcs, we need to alternate the directions in the original path as $(u,v)$, $(w,v)$. This motivates our definition of \textbf{A}-trails and indicates the trend of our next results which show what type of \textbf{A}-trails are conserved by TF-isomorphisms.\\


{\Prop{Let $G$ and $G'$ be mixed graphs and  $P$ be an \emph{\textbf{A}}-trail in $G$. Let $(\alpha,\beta)$ be any  \emph{TF}-isomorphism from $G$ to $G'$.  Then there exists an \emph{\textbf{A}}-trail $P'$ in $G'$ such that $(\alpha,\beta)$ restricted to $P$ maps $P$ to $P'$.}\label{Prop:altertrails01}}\\

{\proof{ For an \textbf{A}-trail consisting of just one arc, the result is trivial. Let us therefore consider an \textbf{A}-trail consisting of $k$ arcs with $k\geq 2$. Let the start vertex of a given \textbf{A}-trail $P
$ be $x_{0}$ and label the successive vertices by $x_{1}, \dots,\ x_{k}$. Assume without loss of generality that $x_{0}$ is the tail of $a_{1}$. The TF-isomorphism maps the arc $a_{1}= (x_{0},x_{1})$ into the arc $a_{1}'=(\alpha(x_{0}),\beta(x_{1}))$. The next arc in $P$ is $a_{2}=(x_{2},x_{1})$ which is mapped by the TF-isomorphism to the arc $a_{2}'=(\alpha(x_{2}),\beta(x_{1}))$ with $\beta(x_{1})$ as a common vertex with $a_{1}'$. By repeating the process until all arcs of $P$ have been included, 
we obtain an \textbf{A}-trail $P'$ of $G'$. Then, by restricting the action of the pair $(\alpha,\beta)$ to $P$ we obtain $P'$ as its image. }}\lesta

This proposition immediately gives the following corollary.\\

{\Cor{If $G$ is a \emph{\textbf{A}}-connected mixed graph which is TF-isomorphic to $H$, then $H$ is also A-connected.}\lesta}

\bigskip

\begin{figure}[h]
 \centering
 \includegraphics[width=8cm,height=8cm]{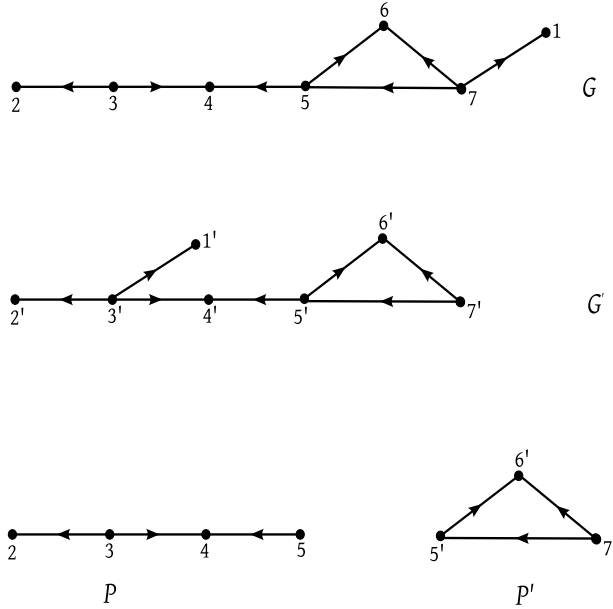}
 \caption{$G$ and $G'$ are TF-isomorphic digraphs but $P$ and $P'$ are not.}\label{fig:nontfisotrails01}
\end{figure}

Proposition \ref{Prop:altertrails01} implies that $Z$-trails are invariant under the action of a TF-isomorphism. The following remarks are aimed to present a clearer 
picture to the reader. Recall that in an \textbf{A}-trail vertices may be repeated so that different alternating trails such as the \textbf{A}-trails $P$ and $P'$  described in Proposition  \ref{Prop:altertrails01}, when taken as digraphs in their own right, may not necessarily be TF-isomorphic. This is illustrated in Figure \ref{fig:nontfisotrails01}.\\

For $G$ and $G'$ as in Figure \ref{fig:nontfisotrails01}, let $\alpha$ map $5$, $7$, $3$ into $5'$, $3'$, $7'$ respectively and let it map arbitrarily the rest of the vertices of $G$ to 
the rest of the vertices of $G'$. Let $\beta$ map $6$, $5$, $1$, $4$, $2$ into $4'$, $2'$, $1'$, $6'$, $5'$ respectively and let it map the rest of the vertices of $G$ to the rest of the 
vertices of $G'$. The maps $\alpha$ and $\beta$ may be represented as shown below where the entries labelled by $*$ may be replaced arbitrarily but without repetitions by any of 
the vertices to which there is no defined mapping.

\[ {\alpha = \left( \begin{tabular}{ccccccc}1&2&3&4&5&6&7\\ \ * &*& 7$'$&*&5$'$&*&3$'$ \end{tabular}                \right)}\qquad
{\beta = \left( \begin{tabular}{ccccccc}1&2&3&4&5&6&7\\1$'$&5$'$&*&6$'$&2$'$&4$'$&* \end{tabular} \right)} \]

\noindent The pair $(\alpha,\beta)$ is then a TF-isomorphism from $G$ to $G'$. However, the alternating trails $P$ and $P'$ in Figure \ref{fig:nontfisotrails01} are not TF-isomorphic digraphs.  On the other hand, as stated in Proposition \ref{Prop:altertrails01}, any \textbf{A}-trail of a given graph, mixed graph or digraph $G$ is mapped by a TF-isomorphism to some \textbf{A}-trail of graph $G'$ whenever $G$ and $G'$ are TF-isomorphic. This is also the case of the trails $P$ and $P'$ in Figure \ref{fig:nontfisotrails01}. In fact it is easy to check  that the open trail $P$ is mapped to the semi-closed trail $P'$ by the pair $(\alpha,\beta)$ as defined above. However $P$ and $P'$ are not TF-isomorphic.\\

{\Prop{Let $G$ and $H$ be mixed graphs. Then a \emph{TF}-isomorphism $(\alpha,\beta)$ from $G$ to $H$ takes closed \textbf{\emph{A}}-trails of $G$ to  closed \textbf{\emph{A}}-trails of $H$.}
\label{prop:tftrailtypefundmap01}}

{\proof{ A closed trail $P$ has an even number of arcs and so it cannot be mapped to a semi-closed trail. Besides, if $P$ were mapped to an open trail,  then $\alpha$ or $\beta$ must map some vertex of $P$ to both the first vertex and last vertex of the open \textbf{A}-trail, which is a contradiction since $\alpha$ and $\beta$ are bijections. As regards the latter case, note that a semi-closed \textbf{A}-trail has an odd number of arcs and a closed \textbf{A}-trail has an even number of arcs.}\lesta}
\bigskip

Therefore, closed \textbf{A}-trails are preserved by TF-isomorphisms just as they are by isomorphisms, but the situation is different for open and semi-closed \textbf{A}-trails.\\

{\Prop{Let $G$ and $H$ be \emph{\textbf{A}}-connected mixed graphs. Then any non-trivial TF-isomorphism $(\alpha,\beta)$ from $G$ to $H$ takes at least one open \textbf{\emph{A}}-trail into a semi-closed \textbf{\emph{A}}-trail and vice-versa.}\label{prop:tftrailtypefundmap}}

{\proof{ As $(\alpha,\beta)$ is non-trivial, there is at least one vertex $u \in$ $V(G)$ such that $\alpha(u)\neq \beta(u)$. Since both $\alpha$ and $\beta$ are bijections, we get $\alpha(u)=\beta(v)$ for some $v\neq u$. Since $G$ is Z-connected, there exists an \textbf{A}-trail joining $u$ and $v$. Clearly $P$ is open.  Its image $P'$ under $(\alpha,\beta)$ is an \textbf{A}-trail of $H$, that starts by $\alpha(u)$ and ends by $\beta(v)$, but since they are equal, $P'$ is semi-closed.\\

Since $(\alpha,\beta)$ is a non-trivial TF-isomorphism from $G$ to $H$, $(\alpha^{-1},\beta^{-1})$ is a non-trivial TF-isomorphism from $H$ to $G$. Therefore, we may use the same arguments to show that $(\alpha^{-1},\beta^{-1})$ must take an open \textbf{A}-trail of $H$ to a semi-closed \textbf{A}-trail of $G$. This implies that $(\alpha,\beta)$ must take some semi-closed \textbf{A}-trail in $G$ to an open \textbf{A}-trail in $H$.}\lesta}

Consider the following example. Let $G$ be a closed \textbf{A}-trail with $6$ vertices and let $H$ consist of a $K_{3}$ and $3$ isolated vertices. Note that $G$ is \textbf{A}-connected whereas $H$ is not. It is straightforward to check that $G$ and $H$ are TF-isomorphic. However, any TF-isomorphism from $G$ to $H$ is clearly non-trivial and maps an open \textbf{A}-trail of length $3$ in $G$ to a semi-closed \textbf{A}-trail of $H$. Therefore, the result of Proposition \ref{prop:tftrailtypefundmap} is false if the hypothesis, namely that both $G$ and $H$ are \textbf{A}-connected, is dropped.\\

As an application of Proposition \ref{prop:tftrailtypefundmap} we get the following result. 

{\Cor{A bipartite graph and a non-bipartite graph cannot be TF-isomorphic. Indeed if $G$ is bipartite and $(\alpha,\beta)$ is a TF-isomorphism from $G$ to some other graph $H$, then $\alpha=\beta$.}\label{cor:tftrailtypefundmap}}

{\proof{Let $G$ be a graph and let $(\alpha,\beta)$ be a non-trivial TF-isomorphism from $G$ to some graph $H$. Then, in view of Proposition \ref{prop:tftrailtypefundmap}, there is an open \textbf{A}-trail of $G$ that is taken to a semi-closed \textbf{A}-trail of $H$. Therefore $H$ has an odd cycle and is non-bipartite. Conversely, $(\alpha^{-1},\beta^{-1})$ is a non-trivial TF-isomorphism from $H$ to $G$ and therefore, by the same argument $G$ cannot be bipartite.}\lesta}

The next section contains a more detailed study of how, using \textbf{A}-trails, a mixed graph $G$ can be made to correspond to a strongly bipartite digraphs, extending the results of Zelinka, particularly those exposed in \cite{zelinka4}. It will turn out that this digraph is a double cover of $G$ which we have already encountered.\\



{\section{Alternating double covers and an equivalence relation on arcs}\label{sec:mzdequival}}

Let $G$ be any mixed graph. Consider the relation $R$ on the set $A(G)$ defined by:
$xRy$ if and only if $x$ and $y$ are the first and last arcs of an \textbf{A}-trail of $G$. Clearly $xRx$ since any given arc is the first
and also the last arc of an \textbf{A}-trail containing only one arc. If $xRy$ then $yRx$ since if $x$ is the first arc of an \textbf{A}-trail,
then $y$ is the last arc and vice-versa. Now suppose that $xRy$ and $yRz$. If $x$ is the first arc of an \textbf{A}-trail $P$, then $y$ is the last arc of
the $P$. Then if $y$ is the first arc of an \textbf{A}-trail $Q$ and $z$ is the last arc of $Q$, then the set-theoretical union of $P$ and $Q$ is an \textbf{A}-trail which has
$x$ as first arc and $z$ as last arc.\\


{\Lem{Let $G$ be a connected graph. Then every two edges of G are joined by trails of both
odd and even length if and only if $G$ is not bipartite.}\label{lem:equivar01}}

{\proof{Let $G$ be non-bipartite. Take any two edges $e_1$, $ e_2$ and fix an odd cycle $C$ of $G$ and two vertices $v_1$, $v_2$ of 
$C$.  Choose two trails $P_1$, $P_2$ joining $e_1$, $e_2$ with $v_1$, $v_2$ respectively. Let $P'$ and $P''$ the two
trails that join $v_1$ and $v_2$ using the edges of $C$. Then there are two trails joining $e_1$, $e_2$, namely $P_1,\ P',\ P_2$ and
$P_1,\  P'',\ P_2$. One of them has odd length and the other even because if $P''$ is odd, $P'$ is even and vice versa and therefore, the inclusion of one instead of the other switches 
parity. Conversely, suppose that there are trails of odd and even length between two fixed edges. In the subgraph induced by these paths, the vertices cannot be partitioned in two 
distinct colour classes and hence this subgraph must be non-bipartite and hence must contain an odd circuit. Hence $G$ contains an odd circuit and is also non-bipartite. }\lesta}

{\Cor{ Let $G$ be a connected graph. Then $R$ has one equivalence class if $G$ is not bipartite and two if it is
bipartite.}\label{Cor:equivalr01}}

{\proof{Suppose that $G$ is  non-bipartite. Consider two arcs $x_1$ and $x_2$ and take the corresponding edges as $e_1$ and $e_2$. Start from $x_1$. On each edge of trails 
joining $e_1$ and $e_2$, choose the arc to obtain an \textbf{A}-trail. Note that before this process may continue for all edges except at most $e_2$. When $e_2$ is reached, the arc corresponding to edge incident with $e_2$ may form an \textbf{A}-trail of order 2 or a directed path. But this depends on whether the concerned trail has odd or even length, so one of the two will give a whole \textbf{A}-trail containing both $x$ and $y$. On the other hand, if $G$ is bipartite, given that $x$ and $y$ form a directed path, which always happens, then every trail joining $x$ and $y$ will be of even length; but an \textbf{A}-trail of even length is open and can't allow directed paths.}\lesta}

\bigskip 

Each equivalence class of $R$ is a set, to which one can naturally associate an \textbf{A}-connected sub-digraph, whose arcs are the elements of the class and whose 
vertices are those incident to at least one of such arcs. In general, the relation $R$ may yield any number of classes, not just one or two as in the case of graphs, as we shall see in Theorem \ref{thm:equivalnewA4} below. \\


If $v$ is any vertex, two different arcs that have $v$ as a first vertex form an \textbf{A}-trail; the same can be said for two different arcs having $v$ as a head. Therefore, the arcs incident with $v$ belong to only one class or two. In the latter case, we say that $v$ is a {\it frontier} vertex. Let $F(G)$ be the set of all frontier vertices of $G$. In view of Corollary \ref{Cor:equivalr01}, if $G$ is a graph then $F(G)$ is either empty (if $G$ is not bipartite) or $F(G)$ $=$ $V$$(G)$ (if $G$ is bipartite). \\

The proof of the next result is straightforward.

{\Prop  {Let $G$ be a connected mixed graph. The following are equivalent:\begin{description}
\item{(i) All classes of $R$ are singletons.}
\item{(ii) All \textbf{\emph{A}}-trails of $G$ are singletons.}
\item{(iii) Each vertex of $G$ has both in-degree and out-degree less than or equal to $1$.}
\item{(iv) $G$ is a directed path or a directed cycle.}\lesta
\end{description}}\label{Prop:equivalnewA1} }

{\Prop{ Let $G$ be a connected mixed graph.  Then $R$ has only one class if and only if the set $F(G)$ is empty. }\label{Prop:equivalnewA2}}

{\proof{The condition is clearly necessary, for if $v$ were an element of $F(G)$ then by definition we would have at least two different classes. On the other hand, if there is more 
than one class, let $x$ and $y$ be arcs that belong to different classes. Since $G$ is connected, there is a trail $P$ that joins $x$ and $y$. Somewhere in $P$ there must be $x'$ and $y'$ that belong to different classes and are incident with a vertex $v$. Thus $v \in F(G)$. }\lesta}

{\Prop{ Let $G$ be a connected mixed graph. Then $F(G)$ is empty or $F(G)$ $=$ $V$$(G)$ or $F(G)$ is a disconnected set of the underlying graph.}\label{Prop:equivalnewA3}}

{\proof{We can assume that $F(G)$ is a proper subset of $V$$(G)$. By Proposition \ref{Prop:equivalnewA2}, there are at least two classes
for $R$. Letting $x$, $y$ be elements of different classes, by the same argument as in Proposition \ref{Prop:equivalnewA2} we infer that
each trial joining $x$ and $y$ must pass through a vertex $v\in F(G)$. Therefore, removing $F(G)$ the arcs $x$ and $y$ end
up in different connected components.}\lesta}

{\Thm{ For every pair $(m,k)$ of positive integers, there exists a mixed graph on which the equivalence relation $R$ induces $m$ classes and having $k$ frontier vertices if and only if $m-1\leq k$.}\label{thm:equivalnewA4}}

{\proof{Let us first construct a mixed graph with $m$ classes and $k$ frontier vertices whenever $m-1 \leq k$. Note that if $m-1=k$ a directed path satisfies the statement (each class consists of a single arc). The same holds for $m-2=k$ and a directed cycle. Assume then that $m-3\leq k$. Consider the $4$-set $\{a,b,c,d\}$ and consider $m-2$ mixed graphs $H_i$ for $i=1,\ \dots,m-2$ where $V$$(H_i)$ $=$  $\{(a,i),(b,i),(c,i),(d,i),(a,i+1)\}$  and $A(H_i)$ contains all the arcs of the triangle $(b,i)$, $(c,i)$, $(d,i)$, plus the additional arcs $((a,i),(b,i))$ and $((d,i),(a,i+1))$. Take any  connected bipartite graph $K$ with $k-m+2$ vertices and fix a vertex $u$ of $K$. Let $L$ be the digraph consisting of the single arc $(u,(a,1))$.\\

Let $G$ be the (standard graph-theoretical) union of $K$, $L$, $H_1,\ \dots, H_{m-2}$. Then $G$ is a connected mixed graph. The classes for $R$ in $G$ are: (i) the class of 
$K$ containing the arcs incident to $u$; (ii) the class of $K$ containing the arcs incident from $u$, together with the extra arc $(u(a,1))$; (iii) each of the $H_i$'s for $i=1,\dots ,\ m-2$. Hence, their number is $m$. Moreover, $F(G)=V(K) \cup \{(a,1),(a,2),...,(a,m-2)\}$, then $|F(G)|=(k-m+2)+(m-2)=k$. Therefore, for all cases where the stated inequality holds, there is a mixed graph $G$ as claimed.\\

Conversely, consider now any mixed graph $G$ and define a graph $X$ such that $V$$(X)$ $=$ $V$$(G)/R$ (that is, the set of classes of $R$ in $G$) and two vertices are adjacent when the associated mixed graphs share a frontier vertex. Then $m=|X|$, while the number $k'$ of edges of $X$ is less or equal to $k=|F(G)|$ (because two classes might share more than a frontier vertex). The known inequality $m-1\leq k'$ implies $m-1 \leq k$ as claimed.}\lesta}

As remarked earlier, a strongly bipartite digraph can be associated with each equivalence class of $R$. Now let these strongly bipartite digraphs  $D_{1}, D_{2}, \dots, D_{k}$ corresponding to the different classes of $R$ obtained from the mixed graph $G$.  Let any vertex $u$ of $V(G)$ which appears as a  source in $D_{i}$ be labelled $u_{0}$ and let any vertex $v$ of $V(G)$ which appears as a sink in $D_{i}$ be labelled $v_{1}$. Therefore, an arc $(u,v)$ in $D_{i}$ now becomes $(u_{0},v_{1})$. It turns out that the strongly bipartite digraph consisting of the components $D_{i}$ labelled this way is \ADC$(G)$ which we have already defined earlier.\\

Figure \ref{fig:mzd04a} shows an example which may be used to illustrate the  following remarks which highlight certain properties of $\mbox{\textbf{ADC}}(G)$ in relation to the mixed graph $G$:\\

\begin{figure}[h]
 \centering
\includegraphics[width = 12.8cm, height = 5.0cm]{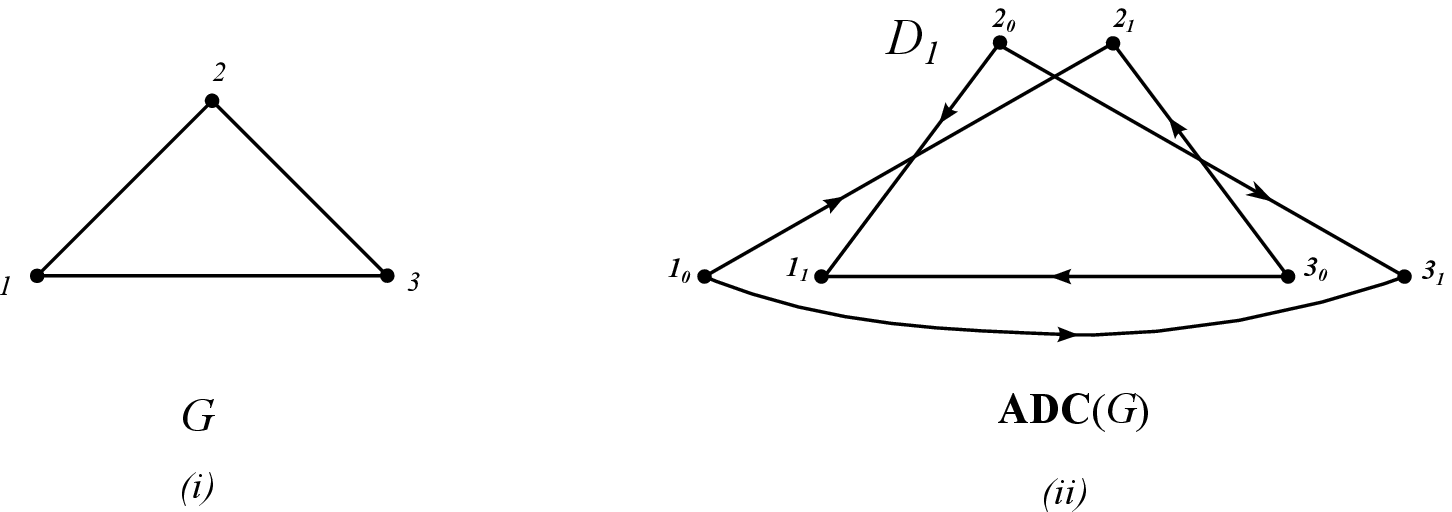}
\caption{$\mbox{\textbf{ADC}}(G)$ obtained form a non-bipartite graph $G$.}\label{fig:mzd01a}
\end{figure}

\begin{figure}[h]
 \centering
\includegraphics[width = 12.8cm, height = 6.4cm]{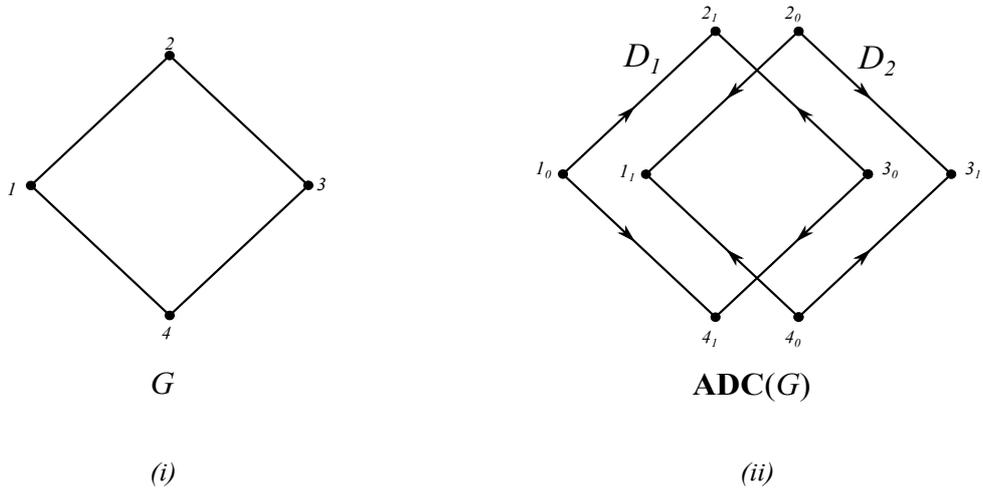}
\caption{$\mbox{\textbf{ADC}}(G)$ obtained from a bipartite graph $G$.}\label{fig:mzd02a} 
\end{figure}

 \begin{description}

\item{1. We know that, for any vertex $u$ of $G$, all incoming arcs $(x,u)$ of $G$ are in the same component of \textbf{ADC}$(G)$ and similarly all outgoing arcs $(u,x)$ of $G$ are in the same component of \textbf{ADC}$(G)$. Therefore $u_0$ if present in \textbf{ADC}$(G)$, cannot appear in two different components. Similarly for $u_1$. However, as we see in examples below, $u_0,$ and $u_1$ can, in some cases, appear in the same component and they can, in other cases, appear in different components. In particular, if $G$ is a bipartite graph  they appear in different components as shown in Figure \ref{fig:mzd02a}$(ii)$ and if $G$ is a non-bipartite graph, they are in the same component as shown in Figure \ref{fig:mzd01a}$(ii)$.}

\item{2. \ADC$(G)$ is a strongly bipartite digraph.}

\item{3. By definition, there is no $u_0$ in \ADC($G$) if $u$ is a sink in $G$, and there is no $u_1$ if $u$ is a source. }
\end{description}

  \medskip

 \begin{figure}[h]
 \centering
\includegraphics[width = 10 cm, height = 11cm]{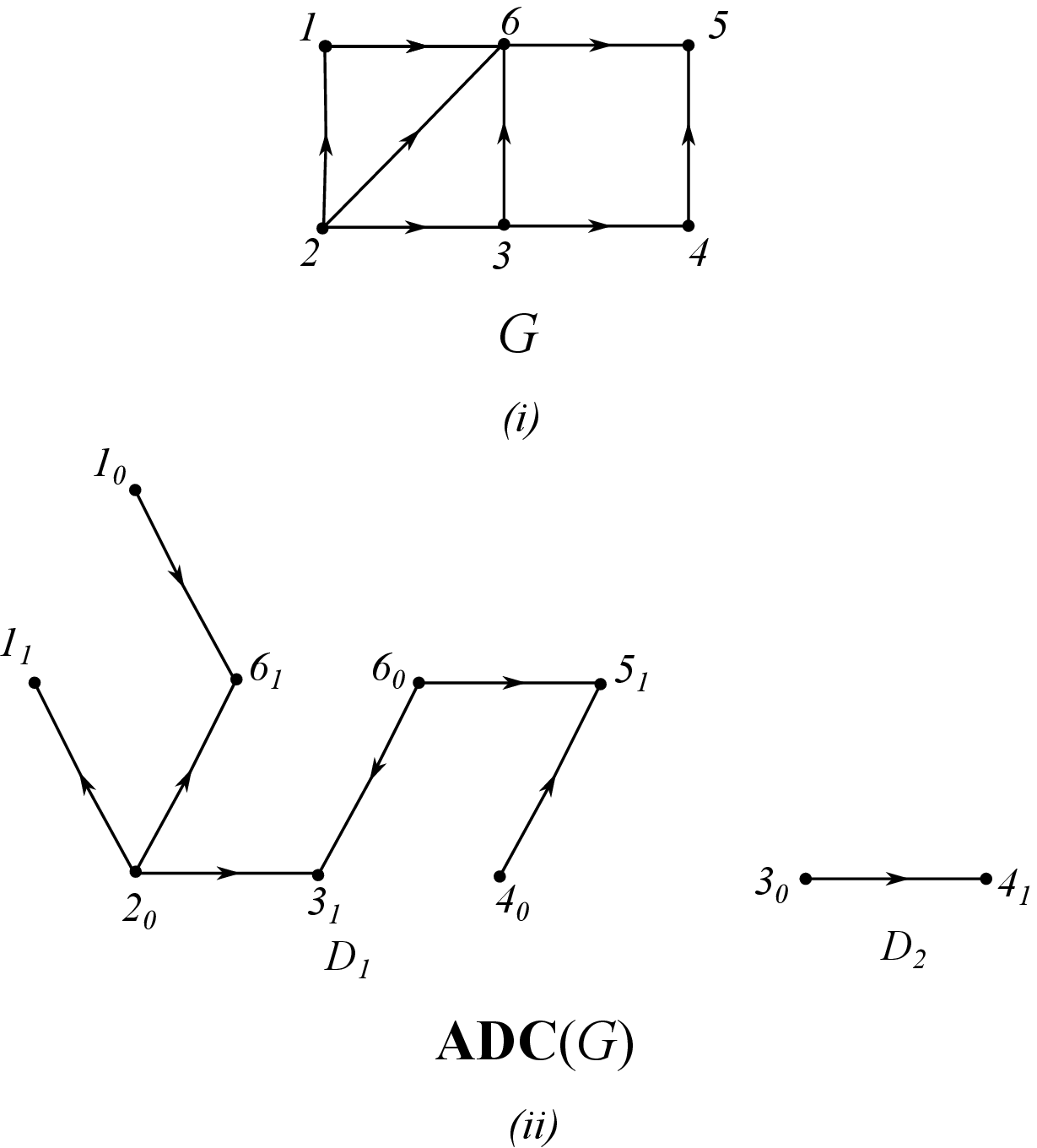}
\caption{$\mbox{\textbf{ADC}}(G)$ obtained from a digraph $G$.}
\end{figure}\label{fig:mzd04a} 
\clearpage

\section{TF-isomorphisms and mixed graph covers}
 

 
The following result can be seen as a corollary to Theorem \ref{thm:idctf} and the proof is easy since the IDC of a mixed graph $G$ can be obtained from $\mbox{\textbf{ADC}}(G)$ simply by removing the directions of the arcs and isolated vertices are irrelevant. Here we give an direct proof because it will help us in later constructions.
 
 {\Thm{Let $G$, $H$ be mixed graphs. The $G$ and $H$ are TF-isomorphic if and only if $\mbox{\textbf{\emph{ADC}}}(G)$ and $\mbox{\textbf{\emph{ADC}}}(H)$ are isomorphic.}\label{Thm:mzdfund01}}
 
{\proof{Let $(\alpha,\beta)$ be a TF-isomorphism from $G$ to $H$. Let $(u,v)$ be an arc of $G$. First note that if $(\alpha(u),\beta(v))$ is an arc of $H$, then $(u_{0},v_{1})$ is an arc of $\mbox{\textbf{ADC}}(G)$ and $(\alpha(u)_{0},\beta(u)_{1})$ is an arc of $\mbox{\textbf{ADC}}(H)$. Let $f$ be a map from $V$$(\mbox{\textbf{ADC}}(G))$ to $V$$(\mbox{\textbf{ADC}}(H))$ such that $f: u_{0} \mapsto x_{0}$ if $x = \alpha(u)$ and $f: v_{1} \mapsto y_{1}$ if $y =\beta(v)$. Consider any arc $(u,v)$ of $G$ and consider the corresponding arc $(u_{0},v_{1})$ in $A(\mbox{\textbf{ADC}}(G))$. Let $(\alpha,\beta)(u,v) = (x,y)$. Then by definition $f$ takes $(u_{0},v_{1})$ to $(x_{0},y_{1})$ in $A(\mbox{\textbf{ADC}}(H))$. The function $f$ maps arcs of $\mbox{\textbf{ADC}}(G)$ to arcs of $\mbox{\textbf{ADC}}(H)$ and it is clearly bijective. Hence, $f$ is an isomorphism from $\mbox{\textbf{ADC}}(G)$ to $\mbox{\textbf{ADC}}(H)$.\\
 
 Now suppose that $\mbox{\textbf{ADC}}(G)$ and $\mbox{\textbf{ADC}}(H)$ are isomorphic. This implies that there exists a map $f$ such that $f(u_{0},v_{1})$ $=$ $(x_{0},y_{1})$. Note that the arcs must always start from a vertex whose label has $0$ as subscript and incident to a vertex whose label has $1$ as subscript, by virtue of the construction presented above. Define $\alpha$, $\beta$ from $V$$(G)$ to $V$$(H)$ as follows. Let $\alpha(u) = x$ if $f(u_{0}) = x_{0}$ where $u \in$ $V$$(G)$ and $x \in$ $V$$(H)$ and let $\beta(v) = y$ if $f(v_{1}) = y_{1}$  where $v \in$ $V$$(G)$ and $y \in$ $V$$(H)$. Then $(\alpha ,\beta)$ takes any arc $(u,v) \in$ $A(G)$ to some $(x,y)$ in $A(H)$. This two-fold mapping is bijective and hence $(\alpha,\beta)$ is a TF-isomorphism from $G$ to $H$.}\lesta}

{\Cor{Let $(\alpha,\beta)$ be a TF-isomorphism from a mixed graph $G$ to a mixed graph $H$. Then there exists an isomorphism  $f_{\alpha,\beta}$ from $\mbox{\textbf{\emph{ADC}}}(G)$ to $\mbox{\textbf{\emph{ADC}}}(H)$ such that $f_{\alpha,\beta}(u_{0},v_{1}) = (x_{0},y_{1})$ if and only if $x = \alpha(u)$ and  $y =\beta(v)$ for some TF-isomorphism $(\alpha,\beta)$ from $G$ to $H$. }\label{cor:zassociated}}

\proof{  The result follows from the proof of Theorem \ref{Thm:mzdfund01}.\lesta}



 Refer to Figure \ref{fig:mzd04}. An isomorphism $f$ from $\mbox{\textbf{ADC}}(G)$ to $\mbox{\textbf{ADC}}(H)$ and the corresponding maps $\alpha$ and $\beta$ from V$(G)$ onto V$(H)$, which are derived from $f$ as described in the proof of Theorem \ref{Thm:mzdfund01}, are given below. \\

{{\begin{eqnarray*}
 f: 1_{0}  \mapsto  1_{0}' &  \ & f: 1_{1}  \mapsto  1_{1}'\\ 
 f: 2_{1}  \mapsto  2_{0}' & \ &  f: 2_{1}  \mapsto  3_{1}'\\ 
 f: 3_{0}  \mapsto 3_{0}' & \ & f: 3_{1} \mapsto  2_{1}'\\ 
 f: 4_{0} \mapsto  6_{0}' &\  & f: 4_{1}  \mapsto  5_{1}'\\ 
 f: 5_{0}  \mapsto  7_{0}' &\  & f: 5_{1}  \mapsto  4_{1}'\\ 
 f: 6_{0}  \mapsto 5_{0}' & \ & f: 6_{1}  \mapsto  6_{1}'\\ 
 f: 7_{0} \mapsto  4_{0}' & \ & f: 7_{1}  \mapsto  7_{1}'
  \end{eqnarray*} }\qquad {{\begin{eqnarray*}
  \alpha: 1  \mapsto  1' &  \ & \beta: 1  \mapsto  1'\\ 
 \alpha: 2  \mapsto  2' & \ &  \beta: 2  \mapsto  3'\\ 
 \alpha: 3  \mapsto 3' & \ & \beta: 3 \mapsto  2'\\ 
 \alpha: 4 \mapsto  6' &\  & \beta: 4 \mapsto  5'\\ 
 \alpha: 5  \mapsto  7' &\  & \beta: 5  \mapsto  4'\\ 
 \alpha: 6 \mapsto 5' & \ & \beta: 6  \mapsto  6'\\ 
 \alpha: 7 \mapsto  4' & \ & \beta: 7  \mapsto  7' 
  \end{eqnarray*} }

 Figure \ref{fig:tod01a} shows a digraph $G$ and its alternating double cover $\mbox{\textbf{ADC}}(G)$ which in this case has three components, namely $D_{1}$, $D_{2}$ and $D_{3}$. Figure \ref{fig:tod01a} also shows how the components of $\mbox{\textbf{ADC}}(G)$ can be combined by associating vertices of the form $u_{0}$ with vertices of the form $v_{1}$, irrespective of whether $u=v$ or $u \not = v$, to form $G$ or other digraphs such as $G_{1}$, $G_{2}$ and $G_{3}$ having the same number of vertices as $G$. It is easy to check that $G$, $G_{1}$, $G_{2}$ and $G_{3}$ are pairwise two-fold isomorphic as expected from the result of Theorem \ref{Thm:mzdfund01} since each of these digraphs have the same number of vertices and have isomorphic \textbf{ADC}s.\\

{\Prop{{(i)} A digraph $H$ is isomorphic to $\mbox{\emph{\textbf{ADC}}}(G)$ for some $G$ if and only if  $H$ is strongly bipartite. 
{(ii)} For every digraph $G$,  ${\mbox{\emph{\textbf{ADC}}}}(\mbox{\emph{\textbf{ADC}}}(G))$ is isomorphic to $\mbox{\emph{\textbf{ADC}}}(G)$.}\label{Prop:gzking01}}
{\proof{We already know that the condition stated in $(i)$ is necessary in order to have $H$ isomorphic to some $\mbox{\textbf{ADC}}(G)$. Conversely, if $H$ has this property, define map $f:$ $V(H)$ $\rightarrow  V (\mbox{\textbf{ADC}}(H))$ as follows: $f(u) = u_{0}$ if $u$ is a source, $f(u) = u_{1}$ if $u$ is a sink. Clearly $f$ is a bijection If $(u,v)$ is an arc of $H$ then by our assumption $u$ is a source and $v$ is a sink of $H$. Then $(u_{0},v_{1})$ $=$ $(f(u),f(v))$ is an arc of $\mbox{\textbf{ADC}}(H)$. Likewise, each arc of $\mbox{\textbf{ADC}}(H)$ takes the form $(u_{0},v_{1})$, with $u$ source and $v$ sink of $H$ and hence $(u_{0},v_{1})$ is the image of $(u,v)$ under $f$. This proves that $f$ is an isomorphism from $H$ to $\mbox{\textbf{ADC}}(H)$, so $(i)$ is satisfied with $G=H$.\\

\noindent Now $(ii)$ is a straightforward consequence of $(i)$, taking $H = \mbox{\textbf{ADC}}(G)$.} \lesta}

 \begin{figure}[h]
 \centering
\includegraphics[width = 14cm, height = 18cm]{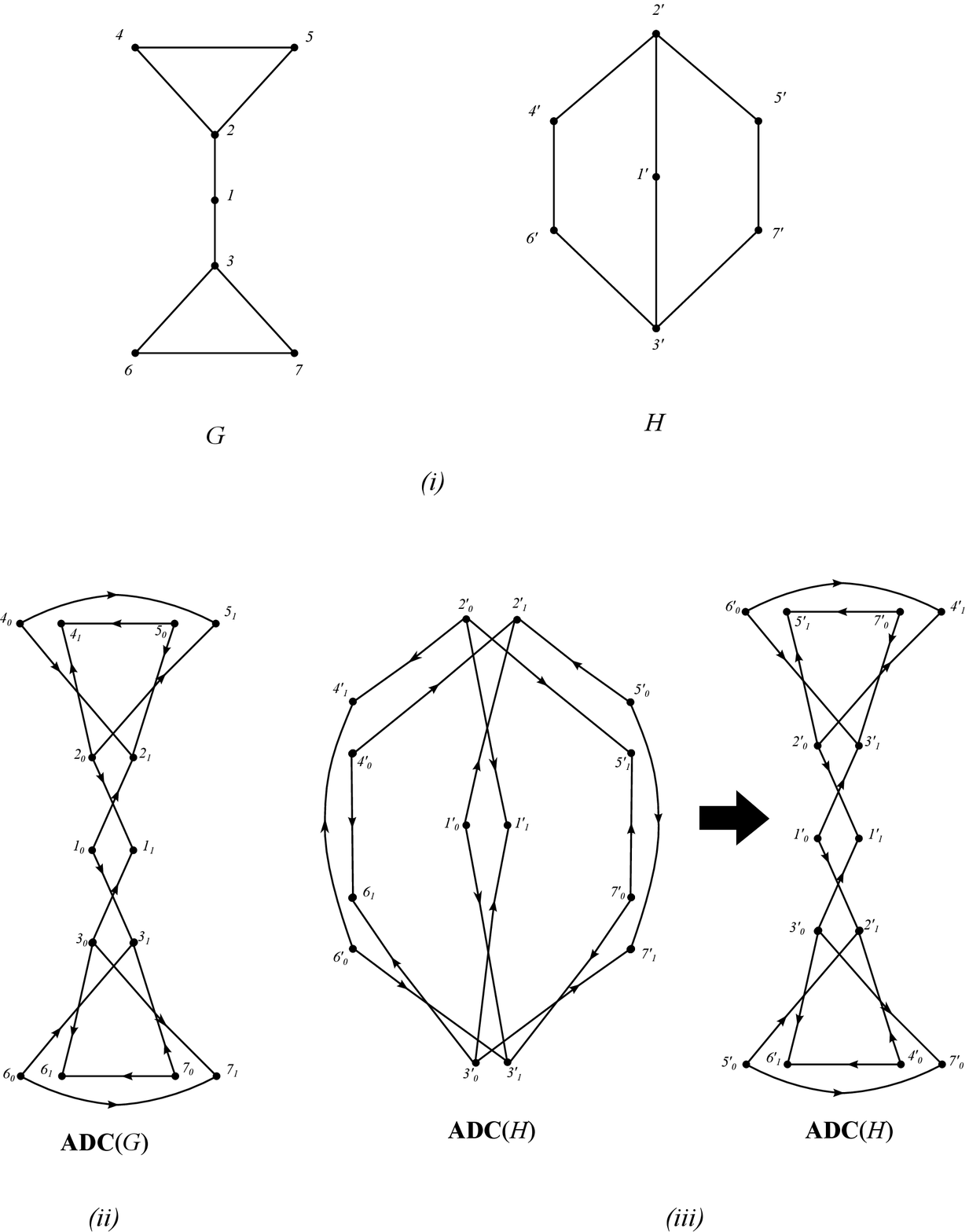}
\caption{$G$ and $H$ are TF-isomorphic graphs and have isomorphic \textbf{ADC}s.}\label{fig:mzd04} 
\end{figure}
\clearpage

 \begin{figure}[h]
 \centering
\includegraphics[width = 14cm, height = 17.5cm]{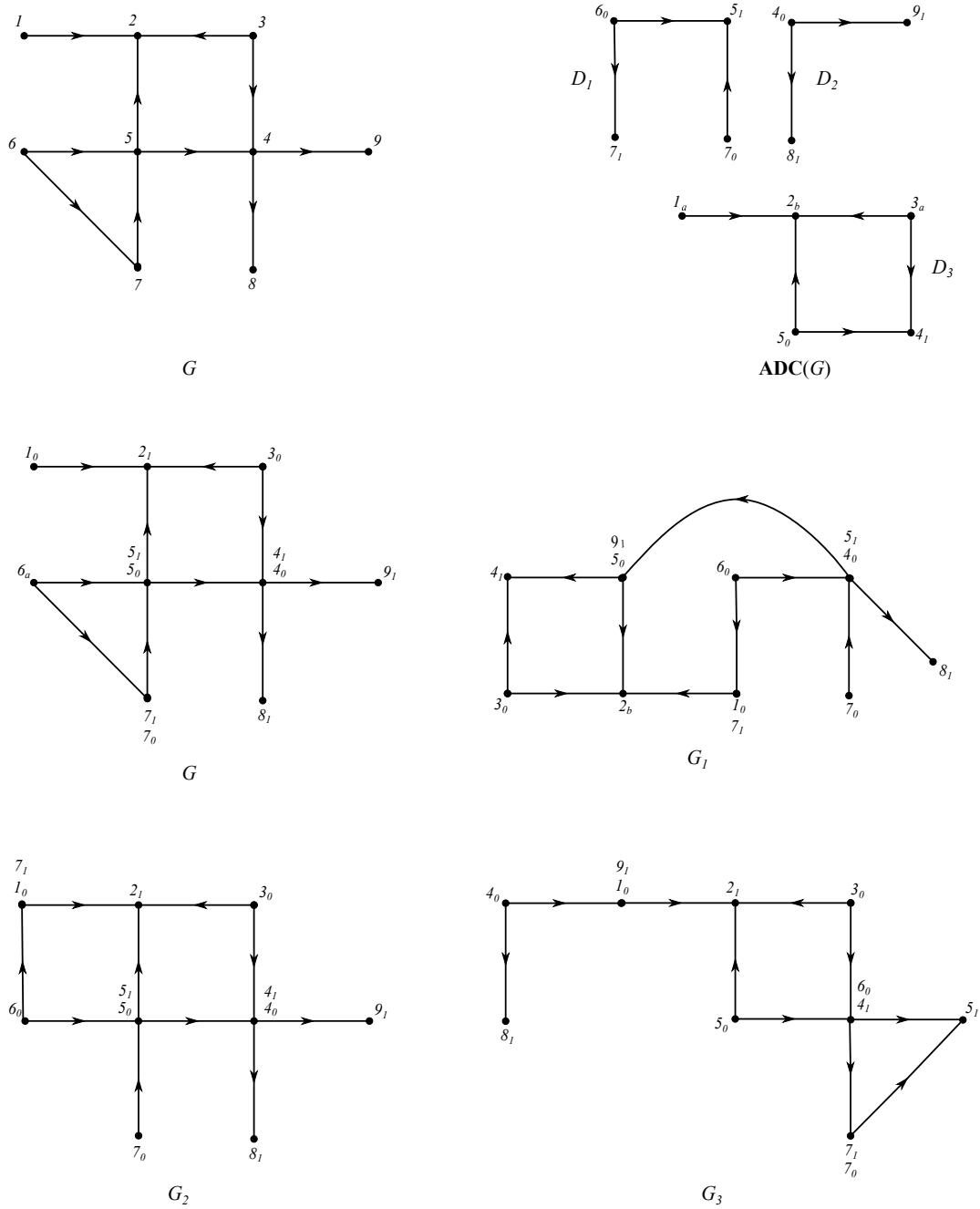}
\caption{$G$ and $H$ are TF-isomorphic graphs and have isomorphic \textbf{ADC}s.}\label{fig:tod01a} 
\end{figure}
\clearpage

\bigskip


\section{Two-fold orbitals}

Let $\mathbf{\Gamma} \leq \mathcal{S} = S_{\mbox{\scriptsize{$|V|$}}} \times  S_{\mbox{\scriptsize{$|V|$}}}$. For a fixed element $(u,v)$ of $V\times V$ let 
\[ \mathbf{\Gamma}(u,v) = \{(\alpha(u),\beta(v)\ |\ (\alpha,\beta) \in \mathbf{\Gamma}\}. \]
The set $\mathbf{\Gamma}(u,v)$ is called a  \emph{two-fold orbital} or TF-\emph{orbital}. A two-fold orbital is the set of arcs of a digraph $G$ having vertex set $V$ which we call \emph{two-fold orbital digraph} or TF-\emph{orbital digraph}. If for every  arc $(x,y)$ in $\mathbf{\Gamma}(u,v)$, the oppositely directed arc $(y,x)$ is also contained in $\mathbf{\Gamma}(u,v)$, then $G$ is a \emph{two-fold orbital graph} or TF-\emph{orbital graph}. This generalisation of the well-known concept of orbital (di)graph has been discussed in \cite{lms1}.\\

 {\Prop{Let $G$ be a strongly bipartite digraph. Then \begin{description}
  \item{(i) There is a homomorphism $\psi$ of \auttfgg onto \emph{Aut}$(G)$.}
  \item{(ii) If $G$ is a TF-orbital digraph, then it is also an orbital digraph.}
  \end{description}}\label{Prop:gzking02}}
  
  {\proof{If $(\alpha,\beta)$ is a TF-automorphism of $G$, define $\psi(\alpha,\beta)$ $=$ $f:$ $V(G) \rightarrow$ $V(G)$ as follows: $f(u) = \alpha(u)$ if $u$ is a source and $f(u)=\beta(u)$ if $u$ is a sink. Since $\alpha$ preserves sources then $f$ takes sources to sources. Similarly, since $\beta$ preserves sinks, then $f$ takes sinks to sinks. Since both $\alpha$ and $\beta$ are permutations, the restrictions of $f$ to the set of sources and to the set of sinks are also permutations. Hence $f$ is a permutation of V$(G)$. Given any arc $(u,v)$ of $G$, note that $(\alpha,\beta)$ takes $(u,v)$ to $(\alpha(u),\beta(v))$, which is equal to $(f(u),f(v))$ because $u$ is a source and $v$ is a sink. Hence $f$ is an automorphism of $G$. so $\psi$ maps \auttfg $\ $ to Aut$(G)$.  A direct computation proves that $\psi$ is a group homomorphism, hence $(i)$ holds.\\
 
Assume now that $G = \mathbf{\Gamma}(u,v)$ for some $\mathbf{\Gamma}$. Then $\mathbf{\Gamma}$ is a subgroup of \auttfg $\ $ and $\psi (\mathbf{\Gamma})$ is a subgroup of Aut $G$. Each arc of $G$ takes the form $(\alpha(u),\beta(v))$, where $(\alpha,\beta)$ $\in$ $\mathbf{\Gamma}$ and $u$, $v$ are a source and a sink respectively. Letting $f = \psi(\alpha,\beta)$ this arc is $(f(u),f(v))$, so it belongs to the orbital digraph $\psi(\mathbf{\Gamma})(u,v)$. This proves that $G$ is contained in this orbital digraph. The opposite inclusion can be shown the same way, so that $G = \psi(\Gamma)(u,v)$ and $(ii)$ follows.\lesta}}\\
  
 {\Cor{Let $G$ be a strongly bipartite digraph. Then $G$ is a two-fold orbital digraph if and only if  ${\emph{\BG}}$ is an orbital digraph.}\label{cor:gzking03}}
 
 {\proof{By Proposition \ref{Prop:gzking01}, $G$ and $\mbox{\textbf{ADC}}(G)$ are isomorphic. If either of them is a TF-orbital, then of course the same holds for the other one, but by Proposition \ref{Prop:gzking02} in this case these TF-orbitals are both orbitals. }\lesta}

 \section{Conclusion}
 
We believe that TF-isomorphisms is a relatively new concept. The only other author who considered them was Zelinka in a short paper motivated by the concept of isotopy in semigroups \cite{zelinka1, Zelinka2}. Our papers (\cite{lms1} and \cite{lms2}) are the first attempts at a systematic study of TF-isomorphisms. \\

In this paper we have shown close links between TF-isomorphisms and double covers, and how the decomposition of a particular double cover can be used to obtain TF-isomorphic graphs.\\

We have also seen that TF-isomorphisms give a new angle for looking at some older problems in graph theory. But does the notion of TF-isomorphism add anything new to these older questions? We believe that it does. For example, in \cite{lms3} we prove this result which explains instability of graphs in terms of TF-automorphisms. \\

{\Thm{Let \emph{Aut}$^{\mbox{\tiny\textbf{{TF}}}}(G)$ be the group of TF-automorphisms of a mixed graph $G$. Then \emph{Aut}$({\emph{\BG}})$ is isomorphic to the semi-direct product \emph{Aut}$^{\mbox{\tiny\textbf{{TF}}}}(G)\rtimes \mathbb{Z}_2$. Therefore $G$ is unstable if and only if it has a non-trivial TF-automorphism.  \lesta}}

Also, it is not very likely that looking at these questions without the notion of TF-isomorphisms would lead one to the notion of \textbf{A}-trails, a technique which we feel is very useful, or the construction of asymmetric graphs with a non-trivial TF-isomorphism, an interesting notion which would be not so natural to formulate using only matrix methods, say. 
Some results and proofs are clearer in the TF-isomorphism setting. For example, in some of the papers cited we find this result about graph reconstruction from neighbourhoods.\\

\begin{Thm}[\cite{Aigner1}] \label{thm:nhoodrecbipartite}
If $G$ is a connected bipartite graph, then any nonisomorphic graph $H$ with the same neighbourhood family as $G$ must be a disconnected graph with two components which themselves have identical neighbourhood hypergraphs. \lesta
\end{Thm}

From the TF-isomorphism point of view, this result follows from three very basic facts: (i) two graphs have the same neighbourhood family (equivalent to being TF-isomorphic) if and only if they have the same canonical double cover; (ii) the canonical double cover of a graph $G$ is disconnected if and only if $G$ is bipartite; and (iii) when $G$ is bipartite, the canonical double cover of $G$ is simply two disjoint copies of $G$. Therefore, for $H$ to have the same canonical double cover as $G$, it must consist of two components isomorphic to $K$, where $G$ is the canonical double cover of $K$. This gives Theorem \ref{thm:nhoodrecbipartite}. And moreover, from these remarks we also see that the only bipartite graphs for which there are non-isomorphic graphs with the same neighbourhood hypergraph are those which are canonical double covers. The Realisability Problem restricted to bipartite graphs therefore becomes:  given a bipartite graph $G$, is there a graph $K$ such that $G$ is the canonical double cover of $K$? A result in this direction was proved in \cite{Scapsalvi2}, where graphs whose canonical double covers are Cayley graphs are characterised\medskip

So it seems that the TF-isomorphism point of view can give a new handle on some of these problems. We intend to pursue this line of research in a forthcoming work.\\

\section*{Acknowledgement}

We are grateful to M.~Muzychuk for first pointing out to us the usefulness of considering TF-isomorphisms as isomorphisms between incidence structures.

 \nocite{porcu}
\nocite{Scapsalvi2}  
  
 \clearpage

\bibliography{reference.bib}
\bibliographystyle{plain}

\end{document}